\documentclass[11pt,svgnames]{amsart}
\usepackage{amsfonts, amssymb, amsmath, amsthm}
\usepackage{mathrsfs} 
\usepackage{latexsym}
\usepackage{enumerate}
\usepackage{multicol}
\usepackage{verbatim}

\usepackage{url}
\usepackage{csquotes}
\usepackage{placeins}
\usepackage{epic}
\usepackage{graphicx}
\usepackage{epstopdf}
\usepackage{pstricks}
\usepackage{pst-plot}
  
\usepackage[colorlinks=true]{hyperref}
\hypersetup{citecolor=blue, linkcolor=blue}
\usepackage{mathabx} 

\input amssym.def
\usepackage{amscd}
\usepackage[mathscr]{eucal}
\usepackage{palatino}
\usepackage{mathtools}
\DeclarePairedDelimiter{\ceil}{\lceil}{\rceil}
\DeclarePairedDelimiter{\floor}{\lfloor}{\rfloor}
\newfont{\cyrr}{wncyr10}
\usepackage{tikz-cd}
\usepackage[utf8]{inputenc}

\newcommand{\thmref}[1]{Theorem~\ref{#1}}

\newcommand{\lemref}[1]{Lemma~\ref{#1}}

\newcommand{\K}{{\mathbf K}}
\newcommand{\F}{{\mathbf F}}

\renewcommand{\b}{{\mathfrak{b}}}
\renewcommand{\a}{{\mathfrak{a}}}
\renewcommand{\c}{{\mathfrak{c}}}

\renewcommand{\P}{{\mathfrak{p}}}

\renewcommand{\O}{{{\mathcal{O}}}}

\newcommand{\rO}{{{\rm{O} }}}
\newcommand{\C}{{\mathbb C}}
\newcommand{\N}{{\mathfrak N}}
\newcommand{\Z}{{\mathbb Z}}
\newcommand{\Q}{{\mathbb Q}}

\renewcommand{\d}{{\partial}}
\newcommand{\q}{{\mathfrak q}}
\newcommand{\R}{{\mathbb R}}

\renewcommand{\L}{\mathcal{L}}

\newcommand{\ogamma}{{\overline{\gamma}}}

\newtheorem{thm}{Theorem}

\newtheorem{lem}[thm]{Lemma}

\newtheorem{rmk}{Remark}[section] 
\newtheorem{defn}{Definition}

\begin{document}
\title{Counting ideals in ray classes}

\author{Sanoli Gun, Olivier Ramar{\'e}  and Jyothsnaa Sivaraman}

\address[Sanoli Gun]   
{Institute of Mathematical Sciences, 
A CI of Homi Bhabha National Institute, 
CIT Campus, Taramani, 
Chennai 600 113, 
India.}

\address[Olivier Ramar{\'e}]
{CNRS / Institut de Math{\'e}matiques de Marseille, 
	Aix Marseille Universit{\'e}, U.M.R. 7373, 
	Site Sud, Campus de Luminy, Case 907, 
	13288 MARSEILLE Cedex 9, France. }

\address[Jyothsnaa Sivaraman]   
{Chennai Mathematical Institute,
	H1, SIPCOT IT Park, Siruseri,
	Kelambakkam,  603103,
	India. }
	
\smallskip

\email{sanoli@imsc.res.in}
\email{olivier.ramare@univ-amu.fr}
\email{jyothsnaas@cmi.ac.in}


\begin{abstract}
Let $\K$ be a number field and $\q$ an integral ideal in $\O_{\K}$.
 A result of Tatuzawa \cite{Tatuzawa*73} from 1973, computes the
 asymptotic (with an error term) for the number of ideals
 with norm at most $x$ in a class of the narrow ray class group
 of $\K$ modulo $\q$. This result bounds the error term with a
 constant whose dependence on $\q$ is explicit but 
 dependence on $\K$ is not explicit.
 The aim of this paper is to prove this asymptotic with a fully 
 explicit bound for the error term.
\end{abstract}

\maketitle   

{\small \tableofcontents}

\section{Introduction and statement of the Theorem}

Given a number field $\K$, the problem of counting the number of
ideals in a given class of the narrow ray class group $H_\q(\K)$
attached to the ideal $\q$ is classical and goes back, if not to Landau,
at least to Hecke. Our query in this paper is the dependence of the
error term on the field $\K$ which we describe fully, and even in a
completely explicit manner.
Let us recall our notation in brief, $n_\K$, $h_{\K, 1}$, $R_\K$, $\mu_\K$ and
$d_\K$ are respectively the degree, the narrow class number, the regulator,
the group of units of finite order in $\K$
and the discriminant of $\K$ while $\alpha_\K$ is the
residue of its Dedekind zeta-function at~$1$. The ring of integers is denoted by $\O_\K$
and $h_{\K}$ denotes the class number of $\K$.

On the technical side, notation $f(x)=\rO^*(g(x))$ means that $|f(x)|\le g(x)$.
In this set-up, we have the following theorem;
\begin{thm}\label{asymfinal}
Let $\q$ be an integral ideal of $\K$ and $[\b]$ be an element of $H_{\q}(\K)$.
For any real number $x \ge 1$, we have
\begin{equation*}
\sum_{\substack{\a\subset\O_\K,\\ [\a]=[\b],\\ \N\a \le x}} 1 = \frac{\alpha_{\K} \varphi(\q)}{|H_{\q}(\K)|}
\frac{x}{\N\q} 
~+~
\rO^* \biggl( E(\K)
  F(\q)^{\frac{1}{n_\K}}\log(3F(\q))^{n_\K}
\left(\frac{x}{\N\q}\right)^{1-\frac{1}{n_\K}}
+ n_{\K}^{8n_\K}\frac{R_\K}{|\mu_\K|} F(\q) \biggr).
\end{equation*}
where $F(\q)=2^{r_1}\varphi(\q)h_\K/h_{\K,\q}$ and
$
  E(\K)=1000
  n_{\K}^{ 12n_{\K}^2 }(R_\K/|\mu_\K|)^{\frac{1}{n_\K}}
  \bigl[\log\bigl((2n_{\K})^{4n_\K}R_\K /|\mu_\K|\bigr)\bigr]^{n_{\K}}$.
\end{thm}
Notice that $F(\q)\ge1$.
Let us briefly recall the definition of the (narrow) ray class group
$H_\q(\K)$.  Let $I(\q)$ be the group of fractional ideals of $\K$
which are co-prime to $\q$ and $P_{\q}$ be the subgroup of $I(\q)$
consisting of principal ideals $(\alpha)$ satisfying
$v_\mathfrak{p}(\alpha -1) \ge v_\mathfrak{p} (\q)$ for all prime
ideals $\mathfrak{p}$ dividing $\q$ and $\sigma(\alpha) >0$ for all
embeddings $\sigma$ of $\K$ in $\R$. We set
$H_{\q}(\K)=I(\q)/P_{\q}$. When $\q=\O_\K$, the group $H_{\q}(\K)$ is
the usual class group in the narrow sense.

The problem of counting ideals in a class of ray class group can be decomposed
in two parts: building a fundamental domain, which turns out to be
made of lattice points in some region, and counting such points.
Our main effort concerns the building of the fundamental domain. Two
hurdles prevents us from directly counting integral ideals of $\K$: the
fact that the narrow ray class group $H_\q(\K)$ is non-trivial, and the
existence of units. To treat both of these, we follow the approach
developed by K.~Debaene in~\cite{Debaene*19}, as it provides us with a
very tame dependence in the field (notice that no discriminant appears
in our error term).  This is combined, as
in~\cite{Debaene*19}, with two general results: the first one shows
that a 'short' enough basis exists for the lattice to be considered,
while the second one counts the lattice points in a given domain.

An overall different approach has
been followed in \cite{Tatuzawa*73} by T.~Tatuzawa, but his results
lack the control of the dependency in $\K$. There also exists
an earlier completely explicit result on this subject, and with a
better error term as far as the dependence in $x$ is concerned. It is
due to J.~Sunley in her PhD memoir \cite{Sunley*71} and is recalled as
Theorem~1.1 of \cite{Sunley*73} (see also \cite{LEE22}). We have few indications as to the
proof of this result, as it has not been published in any journal, but,
knowing that it originates from the method of Landau, we may surmise that
most of the work goes on the dependency of $\K$, relying on a more
classical fundamental domain. We finally mention that the present work
relies on several highly non-trivial results, like the bound
for the regulator given by E.~Friedman in~\cite{Friedman*89}, and the
lower bound for the height of an algebraic number provided by
E.~Dobrowolski in~\cite{Dobrowolski*79}. 

We have several applications of our result which we leave for future
works. The paper is organized as follows. In section \ref{sec2}, we deal with some
notations and preliminaries. In section \ref{sec4}, we
prove \thmref{asymfinal}.

 \section{Notation and Preliminaries}\label{sec2}
\label{Notation}

\subsection*{Notation}
Let $\K\neq \mathbb{Q}$ be a number field with discriminant
$|d_{\K}|\ge3$ (by Minkowski's bound).
Also let $n_{\K} = [\K : \Q]\ge2$ and $\q$ be an (integral) ideal of $\K$. 
The number of real embedding of $\K$ is denoted by $r_1$
whereas the number of complex ones are denoted by $2r_2$. 
The ring of integers of $\K$ is denoted by $\O_{\K}$,
the narrow ray class group modulo $\q$ is denoted by $H_\q(\K)$, its
cardinality by $h_{\K,\q}$
and the (absolute) norm is denoted by $\N$. We shorten $h_{\K,\q}$
by $h_{\K, 1}$ when $\q=\O_\K$.
Whenever required, we shall replace the ideal $\q$ by the
\emph{modulus} $\q_1=\q\q_\infty$, considered as a set of places,
where $\q_\infty$ is the set of all Archimedean places
of~$\K$. But, as we do not consider subsets $S$ of Archimedean places, we may
safely rely only on $\q$ and recall regularly that we count
\emph{narrow} classes. Still to follow tradition, we denote by
$R_{\K,\q_1}$ the $\q_1$-regulator, by $U_{\q_1}$ the corresponding
group of units and by $\mu_{\q_1}$ the number of units of finite order
in $U_{\q_1}$, i.e. also, the cardinality of $\mu_\K\cap U_{\q_1}$.
Throughout the
article $\P$ will denote a prime ideal in $\O_{\K}$ and $p$ will
denote a rational prime number. Further
an element of $H_\q(\K)$ containing an integral ideal
$\a$ will be denoted by $[\a]$.

\subsection*{The Dedekind zeta-function}
For $\Re s= \sigma > 1$, the Dedekind zeta-function is defined by
$$
\zeta_{\K}(s)=\sum_{\a \subseteq \O_{\K} \atop \a \ne 0} \frac{1}{\N(\a)^s}, 
$$
where $\a$ ranges over the integral ideals of $\O_{\K}$. 
It has only a simple pole at $s=1$ of residue $\alpha_{\K}$, say. 
We know from the analytic class number formula that
\begin{equation}\label{acf}
\alpha_{\K} = \frac{2^{r_1} (2\pi)^{r_2} h_{\K} R_{\K}}{|\mu_{\K}| \sqrt{|d_\K|}},
\end{equation}
where $h_{\K}, R_{\K}, d_{\K}$ and $\mu_\K$ are as before. 
\subsection*{The narrow ray-class group}
By narrow ray class group $H_\q(\K)$, we consider that ray class group
where the integral ideal $\q$ is completed with all real
Archimedean places. We have
\begin{equation}\label{eq:4}
|H_1(\K)|\le |H_\q(\K)|\le \varphi(\q)|H_1(\K)|,
\end{equation}
where 
\begin{equation} \label{eq:6}
\varphi(\q)=\N(\q)\prod_{\mathfrak p|\q}\left(1 - \frac{1}{\N(\mathfrak p)}\right)
\end{equation}
and $H_1(\K)$ denotes the narrow ray class group
corresponding to $\O_{\K}$.
A good reference for this are the notes  \cite{Sutherland*15}

by A.~Sutherland.

\subsection*{Orthogonality defect and successive minima}
In this subsection we define some notions and state
results about lattices in $\R^n$ that will be required in due course of the proof.

\begin{defn}
Given a lattice $\Lambda_n$ of rank $n$, the orthogonality
defect $\Omega$ of the lattice $\Lambda_n$ is given by
$$
\Omega = \inf_{(\vec{v_1}, \cdots, \vec{v_n})} 
\frac{||\vec{v_1}|| \cdots ||\vec{v_n}||}{\text{Vol }(\Lambda_n)}
$$ 
where $\{ \vec{v_1}, \cdots, \vec{v_n} \}$ runs over the bases of $\Lambda_n$.
\end{defn}

\begin{defn}
Given a basis $V = \{\vec{v_1}, \cdots ,\vec{v_n}\}$ of a lattice $\Lambda_n$ of rank $n$, 
let $V^{\dagger} = \{\vec{v_1}^{\dagger}, \cdots ,\vec{v_n}^{\dagger}\}$
be the Gram-Schmidt orthogonalisation of $V$.
Let
$$
\alpha_{i,j} = \frac{\vec{v_i} \cdot  \vec{v_j}^{\dagger}}{ || \vec{v_j}^{\dagger}||^2}
\phantom{mm}\text{ for } i, j \in \{1, \cdots,  n \}.
$$
When $n=1$, any basis element of $\Lambda_1$ in $\R$ is defined to be a 
reduced Korkin-Zolotarev basis. When $n>1$, 
the basis $V$ of $\Lambda_n$ is called a reduced Korkin-Zolotarev basis
if it satisfies the following properties;
\begin{enumerate}
\item 
The vector $\vec{v_1}$ is of minimum length among the
vectors $\vec{v_i}$ for $1 \le i \le n$ (with respect to the Euclidean norm),
\item 
The coefficients $|\alpha_{i,1}| \le \frac{1}{2}$  for  $2 \le i \le n$,
\item  
If $\Lambda_{n-1}$ is the orthogonal projection of $\Lambda_n$ 
on the orthogonal complement $(\R\vec{v_1})^{\perp}$, then the 
vectors $\{\vec{v_2} - \alpha_{2,1}\vec{v_1}, \cdots, \vec{v_n} - \alpha_{n,1}\vec{v_1}\}$
also form a  reduced Korkin-Zolotarev basis of $\Lambda_{n-1}$.
\end{enumerate}

\end{defn}
It is easy to see that reduced Korkin-Zolotarev 
bases exist for a lattice $\Lambda_n$ of rank $n$.
\begin{defn}
For a lattice $\Lambda_n$ of rank $n \ge 1$ and for $1 \le i \le n$, the $i$-th successive 
minimum of $\Lambda_n$ is defined by
$$
\delta_i(\Lambda_n) 
= 
\inf\{\lambda \in \R ~|~ B(0, \lambda) \cap \Lambda_n \text{ contains } i \text{ linearly independent vectors}  \}.
$$
Here $B(0, \lambda)$ denotes a ball of radius $\lambda$ around origin in $\R^n$.
For $i=0$, we define $\delta_0(\Lambda_n) =1$.
Further the constant
$$
\gamma_n
= 
\text{sup }\left\{  \left( \frac{\delta_1(\Lambda_n)^n}{  \text{Vol}(\Lambda_n) } \right)^{2/n} ~|~~ 
\Lambda_n \text{ is a lattice of rank } n  \right\}
$$
is called the Hermite's constant.
\end{defn}
In this set-up, we have the following theorem.
\begin{thm} [Lagarias, Lenstra and Schnorr \cite{Lagarias-Lenstra-Schnorr*90}]\label{KZ}
If $\{\vec{v_1} \ldots \vec{v_n}\}$ is a reduced Korkin-Zolotarev basis
for a lattice $\Lambda_n$ of rank $n$, then
$$
\prod_{j=1}^{n} ||\vec{v_j}||^2 
\le \left(\prod_{j=1}^n \frac{j+3}{4} \right) \gamma_n^n \text{ Vol}(\Lambda_n)^2,
$$ 
where $\gamma_n$ is the Hermite's constant.
Further, an upper bound for the Hermite's constant is given by
$$
\gamma_n 
\le
n  \phantom{mm}\text{ for } n \ge 1.
$$
\end{thm}

We now define the notion of the Lipschitz class of a subset of $\R^n$.
\begin{defn}\label{def-2}
Let $S$ be a subset of $\R^n$ with $n \ge 2$. We
say that $S$ is of Lipschitz class $\L(n,M,L)$
if there are $M$ maps $\phi_1,\ldots \phi_M : [0,1]^{n-1} \to \R^n$
such that $S$ is contained in the union of images of $\phi_i$
for $i\in \{1,\ldots M\}$ and
$$
||\phi_i(\overline{x}) - \phi_i( \overline{y})|| \le L  ~||\overline{x} -  \overline{y}||,
$$
where $\overline{x}, \overline{y} \in [0,1]^{n-1}$.
\end{defn}

We conclude by stating a theorem of Widmer \cite{Widmer*10a} which allows us to 
estimate the main term as well as the error term in \thmref{counting}.

\begin{thm}[Widmer \cite{Widmer*10a}]\label{Widmer}
Let $\Lambda_n$ be a lattice in $\R^n$ with successive minima 
$\delta_0(\Lambda_n),  \cdots,  \delta_n (\Lambda_n)$. 
Let $S$ be a bounded set in $\R^n$ such that its boundary
is of Lipschitz class $\L(n,M,L)$ for some natural number $M$ and
positive constant $L$.  Then $S$ is measurable and 
$$
\left|  |S \cap \Lambda_n| - \frac{\text{Vol}(S)}{\text{Vol}(\Lambda_n)} \right| \le Mn^{3n^2/2} 
\max_{0 \le i < n} \frac{L^i}{\delta_0(\Lambda_n) \cdots \delta_i(\Lambda_n)}.  
$$
\end{thm}
\subsection*{Lower bounds for algebraic conjugates}
In this subsection we recall a theorem of Dobrowolski which
gives a lower bound on the absolute value of all the conjugates
of an algebraic integer which is not zero or a root of unity.
\begin{thm} [Dobrowolski \cite{Dobrowolski*79}\label{Dob}]
Let $\alpha$ be a non-zero algebraic integer of degree $n > 1$ 
and let $\tilde{\alpha}$
be the maximum of the absolute values of all conjugates of $\alpha$.
If $\alpha$ is not a root of unity, then
$$
\tilde{\alpha}
\ge 
1 + \frac{\log n}{6n^2}.
$$
\end{thm}

\smallskip

\section{Counting integral ideals in classes of the ray class group}\label{sec4}

Let $\{\sigma_1, \cdots \sigma_{n_{\K}} \}$ be the set of all embeddings 
of $\K$ into $\C$.  The first $r_1$ embeddings are all real embeddings 
and the embeddings $\{\sigma_{r_1 + i},  ~ \sigma_{r_1 + r_2 + i} \}$
for $1 \le i \le r_2$ are complex conjugates.
Consider the first $r_1 + r_2$ embeddings from
this set. We will use $r$ to denote $r_1+r_2-1$ and as before
$\q_1 = \q \q_{\infty}$ to denote a modulus, where $\q \subseteq \O_{\K}$ be an ideal
and $\q_{\infty}$ contains all real places of $\K$.

\medskip
\subsection{Fundamental domain}
Let $\O_{\K}^*$ be the group of units of $\O_{\K}$ and $U_{\q}$ (respectively $U_{\q_1}$) 
be the subgroup of $\O_{\K}^*$ consisting of units which are $1 \bmod \q$ 
(respectively $1 \bmod^* \q_1$).  Both of these subgroups 
$U_{\q}$ and $U_{\q_1}$ are of finite index in $\O_{\K}^*$. 
Let $\phi$ denote the embedding
\begin{eqnarray*}
\phi  :   \K  & \to & \R^{r_1} \times \C^{r_2} \\
 x  & \to & ( \sigma_i(x) )_{i=1}^{r +1}
\end{eqnarray*} 
Further let $f$ denote the map
\begin{eqnarray*}
f : \R^{r_1} \times \C^{r_2} & \to & \R^{r+1} \\
(x_i)_i & \to & (\log |x_i|)_{i=1}^{r + 1}
\end{eqnarray*}
Since $[\O_{\K}^* : U_{\q_1}]$ is finite, the image under the map
$f \circ \phi$ of $U_{\q_1}$ is also a lattice of rank $r$. Let
$\{\eta_1, \ldots \eta_{r}\}$ be a set of multiplicatively independent
generators for the group $U_{\q_1}$ modulo roots of unity and the vectors
\begin{align*}
\vec{v_1} 
& = 
 \left(\frac{1}{n_{\K}}, \cdots ,\frac{1}{n_{\K}}\right), 
&
\vec{v_j} 
 = 
\left( \log |\sigma_i(\eta_{j-1})|\right)_{i=1}^{r+1}
\phantom{m} \text{for  } 2 \le j \le r+1
\end{align*}
form a basis for $\R^{r+1}$.  The vectors $\vec{v_2}, \cdots \vec{v_{r+1}}$
form a basis for a lattice of rank $r$ and 
$\vec{v_1}, \vec{v_2}, \cdots \vec{v_{r+1}}$ are $\R$ linearly independent. 
We can now write the vector 
$( \log |x_i| )_{i=1}^{r +1}$ as
$$
( \log |x_i|)_{i=1}^{r +1} = \alpha_1(x) \vec{v_1} 
+ \cdots + \alpha_{r+1}(x) \vec{v_{r+1}},
$$
and therefore
\begin{equation}\label{id}
|x_i| = e^{\alpha_1(x)/n_{\K}} \prod_{j=2}^{r+1} 
|\sigma_i(\eta_{j-1})|^{\alpha_j(x)} 
= \alpha(x)^{1/n_{\K}} \prod_{j=2}^{r+1} 
|\sigma_i(\eta_{j-1})|^{\alpha_j(x)},
\end{equation}
where $x = (x_1, \cdots, x_{r+1})$ and $\alpha(x) = e^{\alpha_1(x)}$.
If $x \in \K$, then $x_i = \sigma_i(x)$ for $1 \le i \le r+1$ and so taking product over
all $i$ in \eqref{id}, we identify $\alpha(x)$:
\begin{equation}\label{id-K}
 \prod_{i=1}^{r +1 } |\sigma_i(x)|^{e_i}  = |\N(x)| =\alpha(x),
\end{equation}
where $e_i=1$ when $i\le r_1$ and $e_i=2$ otherwise.
We now define the map
\begin{eqnarray*}
g : \R^{r_1} \times \C^{r_2} & \to & \R^{r+1}\\
x= (x_i)_i & \to & (\alpha(x), \alpha_2 (x), \cdots \alpha_{r+1}(x)).
\end{eqnarray*}
We now define $\F = g^{-1}(\R_{+} \times [0,1)^{r})$. 
This corresponds to the set $(f \circ \phi)^{-1} (S')$,
where $S'$ is given by the points corresponding to the vectors
$$
\{  \alpha_1\vec{v_1} + \ldots + \alpha_{r+1} \vec{v_{r+1}} 
~\mid~ \alpha_1 \in \R, ~\alpha_i \in [0,1) \text{ for  } i >1\}.
$$
Since the vectors give rise to a lattice of full rank,
given an $x \in \K$, there is an
$\eta \in U_{\q_1}$ such that $\phi(x/\eta) \in \F$.
Conversely, for an $x \in \K$,
with $\phi(x) \in \F$ and any $\eta \in U_{\q_1}$, we note that 
$\phi(\eta  x) \in \F$ if and only if $\eta$ is a root of unity.

\subsection{Notation}
Throughout the rest of the section, we will use the following notations. 
\begin{enumerate}
\item $\F(a_1, b_1, \ldots, a_r, b_r, X) = g^{-1}\left( (0,X] \times \prod_{j=1}^{r} [a_j, b_j)\right)$.

\item $\F_{\frac{1}{2}}(a_1, b_1, \ldots, a_r, b_r, X) = g^{-1}\left( \left(\frac{X}{2},X\right] \times
 \prod_{j=1}^{r} [a_j, b_j)\right)$. 

\item For $\overline{\gamma} = ( \gamma_i)_{i=1}^{r_1} \in \{\pm 1\}^{r_1}$, 
denote by $\R^{r_1}_{ \bar{\gamma}} = \{ (x_1, \cdots, x_{r_1}) ~|~ \text{sign} (x_i) = \gamma_i\}$.
Then 
$\F_{\overline{\gamma}}(a_1, b_1, \ldots, a_r, b_r, X) 
= \F(a_1, b_1, \ldots, a_r, b_r, X) \cap (\R_{\overline{\gamma}}^{r_1} \times \C^{r_2})$. \\
Further,  $\F_{\frac{1}{2},\overline{\gamma}}(a_1, b_1, \ldots, a_r, b_r, X) 
= \F_{\frac{1}{2}}(a_1, b_1, \cdots, a_r, b_r, X) \cap (\R_{\overline{\gamma}}^{r_1} \times \C^{r_2})$.

\item 
$\F_{\frac{1}{2}, \overline{\gamma}}(X) = \F_{\frac{1}{2}, \overline{\gamma}}(0,1,\cdots 0,1,X)$.

\item For $j \in \{1, \ldots, r\}$, $m_j = \max_i \ceil{\log |\sigma_i(\eta_j)|}$,
where $\{\eta_1, \ldots \eta_{r}\}$ is a set of multiplicatively independent
generators for the group $U_{\q_1}$ modulo roots of unity.
\end{enumerate}

\medskip
\subsection{Computing the Lipschitz class of the boundary} \label{lipschitzclass}
\begin{defn}
Let $\q_1$ be a modulus of $\K$ and $U_{\q_1}$ be the subgroup
of units which are $1 \bmod^* \q_1$. Further let $\{\sigma_1, \cdots, \sigma_{r+1}\}$
be a set of $r+1$ embeddings of $\K$ into $\C$ 
where no two embeddings are conjugates of each other.
If $\{ \eta_1 , \cdots, \eta_r\}$ are multiplicatively independent
units generating $U_{\q_1}$ modulo the roots of unity, then the $\q_1$ regulator
$R_{\K,\q_1}$ is defined by
$$
R_{\K,\q_1} = \left| {\rm det} \left(e_i\log|\sigma_i(\eta_j )|  \right)_{i,j} \right|
$$
where $i,j \in \{1,2 \ldots r\}$ and $e_i= 1$ or $2$ if the corresponding 
embedding is real or complex respectively. 
\end{defn}
Since $\sum_{i=1}^{r+1} e_i\log|\sigma_i(\eta_j)| =0$
for $1 \le j\le r$, we see that $R_{\K,\q_1}$ is independent of the
generating set of units and the choice of $r$ embeddings.

\smallskip

\begin{lem}\label{SPT}
Let $\Gamma$ be a set of points in $\R^{r_1} \times \C^{r_2}$,
$\vec{k} = (k_j)_{j=1}^{r} \in \prod_{j=1}^{r} ([0,m_j)~\cap~\Z)$ and
$$
\beta_{\vec{k}} = \left( \prod_{j=1}^r |\sigma_i(\eta_j)|^{-k_j/m_j}\right)_{i=1}^{r+1}.
$$
Then for any $\overline{\gamma} = ( \gamma_i)_{i=1}^{r_1} \in \{\pm 1\}^{r_1}$
and any positive real $X$, we have
$$
|\Gamma \cap \F_{\frac{1}{2},\overline{\gamma}}(X)| 
= 
\sum_{\vec{k}} \bigg|(\Gamma \cdot \beta_{\vec{k}}) 
\cap 
\F_{\frac{1}{2},\overline{\gamma}}\left(0, \frac{1}{m_1}, \cdots 0, \frac{1}{m_r}, X\right)\bigg|,
$$
where the sum runs over $\vec{k} \in \prod_{j=1}^{r} ([0,m_j)~\cap~\Z)$.
\end{lem}

\begin{proof}
We have
\begin{equation*}
|\Gamma \cap \F_{\frac{1}{2},\overline{\gamma}}(X)| 
 = 
\sum_{(k_1,\cdots, k_r) \in \Z^r, \atop{0 \le k_i \le m_i-1}} 
\bigg|\Gamma \cap \F_{\frac{1}{2},\overline{\gamma}}\left(\frac{k_1}{m_1}, \frac{k_1+1}{m_1}, 
\cdots \frac{k_r}{m_r}, \frac{k_r+1}{m_r}, X\right)\bigg|.
\end{equation*}
If we multiply an element of $\F_{\frac{1}{2}, \overline{\gamma}} \left(\frac{k_1}{m_1}, 
\frac{k_1+1}{m_1}, \cdots \frac{k_r}{m_r}, \frac{k_r+1}{m_r}, X\right)$ 
by $\beta_{\vec{k}= (k_1,\ldots k_r)}$, using~\eqref{id},
we see that we get an element in $\F_{\frac{1}{2}, \overline{\gamma}} 
\left(0, \frac{1}{m_1}, \cdots 0, \frac{1}{m_r}, X\right)$. 
Therefore
\begin{equation*}
|\Gamma \cap \F_{\frac{1}{2},\overline{\gamma}}(X)| 
 = 
\sum_{(k_1,\ldots, k_r) \in \Z^r, \atop{0 \le k_i \le m_i-1}} 
\bigg| (\Gamma \cdot \beta_{\vec{k}}) \cap \F_{\frac{1}{2},\overline{\gamma}}
\left(0, \frac{1}{m_1}, \ldots 0, \frac{1}{m_r}, X\right)\bigg|.
\end{equation*}
\end{proof}
Let $h' : \R^{r_1} \times \C^{r_2} \to \R^{n_{\K}}$ be a map defined by
$h' (y_1, \cdots,  y_{r_1+r_2}) = (z_1, \cdots, z_{n_{\K}} )$, 
where for $1 \le i \le r_1$, $z_i = y_i$ and for $r_1+1 \le i \le r_1 + r_2$, 
$$
z_i = |y_i|\cos (\arg(y_i))
\phantom{m}\text{and}\phantom{m}
z_{r_2+i} =  |y_i| \sin (\arg(y_i)),
$$
where argument for $y_i$'s are inside $[0, 2\pi)$.  

\begin{lem}\label{length}
Let $X$ be any positive real number and 
$$
x = (x_1, \cdots x_{r_1+r_2}) \in \F \left(0,\frac{1}{m_1}, \cdots, 0, \frac{1}{m_r}, X \right).
$$
Then $|| h' (x) || ~\le~ (\sqrt{r+1})~e^r~ X^{1/n_{\K}}$,
where $||\cdot||$ is the Euclidean norm on $\R^{n_{\K}}$.
\end{lem}
\begin{proof}
From the definition of the norm, we have
$$
|| h' (x) || 
=
 \sqrt{\sum_{i=1}^{r_1} |x_i|^2 + \sum_{i= r_1+1}^{r_1+r_2} (|\Re(x_i)|^2 + 
 |\Im(x_i)|^2)}.
$$
From \eqref{id} and the definition of $\F(0,\frac{1}{m_1}, \cdots, 0, \frac{1}{m_r}, X)$, 
we see that
$$
|x_i|  = |\alpha(x)|^{1 / n_{\K}} \prod_{k=2}^{r+1} 
|\sigma_i(\eta_{k-1})|^{\alpha_k(x)}
$$
where $0 < \alpha(x)  \le  X$ and $0 \le \alpha_k(x) < 1/m_k$ for $2 \le k \le r + 1$. 
Hence for $1 \le i \le r+1$, we have 
$$
|| h' (x) || ~\le~ (\sqrt{r+1}~)~e^r~ X^{1/n_{\K}} .
$$
This completes the proof of the lemma.
\end{proof}
The next lemma can be found in \cite[Lemma 6]{Debaene*19}.
\begin{lem}\label{lipschitz}
Let $f : \R^m \to \R$ be a function such that $f(y_1,\ldots y_m) = c \prod_{j=1}^m g_j(y_j)$,
where $g_j : \R \to \R$ are functions that satisfy $|g_j(y_j) - g_j(y_j')| \le K_j|y_j - y_j'|$ and
$|g_j(y_j)| \le M_j$. Then we have
$$
|f(\overline{y}) - f(\overline{y}')| 
\le 
\left( c\sum_{j=1}^m K_j \prod_{k \neq j} M_k\right) ~ ||\overline{y} - \overline{y}'||,
$$
where $\overline{y}, \overline{y}' \in \R^m$ and $||.||$ is the Euclidean norm on $\R^m$.
Let $h : [0,1]^{m-1} \to \R^m$, $h=(h_1, \ldots h_m)$ where $h_i : [0,1]^{m-1} \to \R$
are functions that satisfy 
$|h_i(\overline{y}) - h_i(\overline{y}')| \le L_i ~||\overline{y} - \overline{y}'||$. 
Then we have 
$$
|h(\overline{y}) - h(\overline{y}')| \le \sqrt{m}  ~(\text{max}_i L_i )~ ||\overline{y} - \overline{y}'||,
$$
where $\overline{y}, \overline{y}' \in [0,1]^{m-1}$.
\end{lem}
We now compute the Lipschitz class of the boundary of the fundamental domain $\F$.
From now on, we will denote the boundary of a set $\F$ by $\d\F$.
\begin{lem}\label{Lip-const}
For a number field $\K \neq \Q$, $\overline{\gamma}  \in \{\pm 1\}^{r_1}$ and positive 
real number $t$, the set 
$$
\d\F_{\frac{1}{2}, \overline{\gamma}}\left(0,\frac{1}{m_1}, \cdots, 0, \frac{1}{m_r}, t^{r+1} \right)
$$
is of Lipschitz class $\L(n_{\K}, 2r+2,  Lt)$, where
$L = \sqrt{n_{\K}}(2\pi+r)e^r$.
\end{lem}
\begin{proof}
Let $x= (x_i)_{i=1}^{r+1} \in \F_{\frac{1}{2}}(0,\frac{1}{m_1}, \cdots, 0,  \frac{1}{m_r}, 1)$.
Using \eqref{id}, we see that $|x_i|$ for all $i$ is 
uniformly bounded away from $0$ by a constant 
that depends only on $\K$ and not on $i$.  Hence
$$
\F_{\frac{1}{2}} \left(0,\frac{1}{m_1}, \cdots, 0,  \frac{1}{m_r}, 1\right) 
~\subseteq~
\bigcup_{\overline{\gamma}}(\R_{\overline{\gamma}}^{r_1} \times \C^{r_2}),
$$
where $\overline{\gamma}$ varies over elements of $\{ \pm1 \}^{r_1}$.
Since $(\R_{\overline{\gamma}}^{r_1} \times \C^{r_2})$ are disjoint for
distinct $\overline{\gamma}$, we have
$$
\F_{\frac{1}{2}}\left(0,\frac{1}{m_1}, \cdots, 0,  \frac{1}{m_r}, 1\right) 
=
\bigcup_{\overline{\gamma}}\F_{\frac{1}{2}, \overline{\gamma}}
\left(0,\frac{1}{m_1}, \cdots, 0,  \frac{1}{m_r}, 1\right).
$$
Since $\overline{\F_{\frac{1}{2}}\left(0,\frac{1}{m_1}, \cdots, 0,  \frac{1}{m_r}, 1\right)}$
does not contain any point $x= (x_i)_{i}$ with $x_i=0$ for some $i$, we
note that $\overline{\F_{\frac{1}{2}, \overline{\gamma}}
\left(0,\frac{1}{m_1}, \cdots, 0,  \frac{1}{m_r}, 1\right)}$ are
disjoint for distinct $\overline{\gamma}$. This implies that

\begin{equation}\label{f-id}
\d\F_{\frac{1}{2}}\left(0,\frac{1}{m_1}, \cdots, 0,  \frac{1}{m_r}, 1\right) 
=
\bigcup_{\overline{\gamma}}\d\F_{\frac{1}{2}, \overline{\gamma}}
\left(0,\frac{1}{m_1}, \cdots, 0,  \frac{1}{m_r}, 1\right).
\end{equation}
More precisely
$$
\d\F_{\frac{1}{2}, \overline{\gamma}} \left(0,\frac{1}{m_1}, \cdots, 0, \frac{1}{m_r}, 1\right)
 = \d\F_{\frac{1}{2}}\left(0,\frac{1}{m_1}, \cdots , 0, \frac{1}{m_r}, 1\right)
\cap (\R_{\overline{\gamma}}^{r_1} \times \C^{r_2}).
$$
We want to compute the Lipschitz constant for the set
$\d\F_{\frac{1}{2}, \overline{\gamma}}(0,\frac{1}{m_1}, \cdots, 0, \frac{1}{m_r}, 1)$.
For any element $\overline{\gamma} = (\gamma_i)_{i=1}^{r_1}\in \{\pm 1\}^{r_1}$, consider the map
\begin{eqnarray*}
\tilde{g}_{\ogamma}  :  \R_{\overline{\gamma}}^{r_1}  \times  \C^{r_2}  & \to & \R^{n_{\K}}\\
 (y_1, \cdots, y_{r_1+r_2}) & \to & \left( g(y_1,\cdots y_{r_1+r_2}), \arg(y_{r_1+1}), \cdots, \arg(y_{r_1+r_2})\right)
\end{eqnarray*}
where the argument is chosen from $[0,2\pi)$.
Then
\begin{align*}
\d\F_{\frac{1}{2}}\left(0,\frac{1}{m_1}, \cdots, 0, \frac{1}{m_r}, 1\right) 
\cap (\R_{\overline{\gamma}}^{r_1} \times \C^{r_2})
& = 
\d g^{-1} \left( [\frac{1}{2}, 1) \times \prod_{i=1}^{r} [0,\frac{1}{m_i}) \right) 
                                                      \cap (\R_{\overline{\gamma}}^{r_1} \times \C^{r_2})
  \\
& = 
\d {\tilde{g}_{\ogamma}}^{-1} \left( [\frac{1}{2},1) 
\times \prod_{i=1}^{r} [0,\frac{1}{m_i}) \times [0,2\pi)^{r_2} \right).
\end{align*}
The map
\begin{eqnarray*}
\tilde{g}_{\ogamma} : \R_{\overline{\gamma}}^{r_1}  \times  \C^{r_2} 
& \to & 
\R^{r+1} \times [0, 2\pi)^{r_2} \\
(y_1, \cdots, y_{r_1+r_2}) 
& \to & 
\left(g(y_1, \cdots, y_{r_1 + r_2} ), \arg(y_{r_1+1}), \cdots, \arg(y_{r_1+r_2})\right)
\end{eqnarray*}
is a bijection. Indeed if 
$\tilde{g}_{\ogamma}(y_1, \ldots, y_{r_1+r_2}) = \tilde{g}_{\ogamma}(y'_1, \ldots, y'_{r_1+r_2})$,
then $|y_i| = |y'_i|$ for $1 \le i \le r_1+r_2$. Further the first $r_1$ real numbers have the same sign
and the last $r_2$ complex numbers have the same arguments. Hence
$(y_1, \ldots, y_{r_1+r_2} ) = (y'_1, \ldots, y'_{r_1+r_2})$.
When  $(x_1, \ldots x_{n_{\K}}) \in R^{n_{\K}}$, then
for $1 \le i \le r_1$, define
\begin{eqnarray} \label{inv1}
y_i = \gamma_i x_1^{1/n_{\K}} \prod_{k=2}^{r+1} |\sigma_i(\eta_{k-1})|^{x_k}
\end{eqnarray}
and for $ r_1 + 1 \le i \le r_1+r_2$, define
\begin{equation} \label{inv2}
y_i = \left(x_1^{1/n_{\K}} \prod_{k=2}^{r+1} |\sigma_i(\eta_{k-1})|^{x_k}\right) 
e^{2i\pi x_{r_2 + i}} .
\end{equation}
Clearly $(y_1, \ldots, y_{r_1+r_2} ) \in \R_{\overline{\gamma}}^{r_1}  \times  \C^{r_2} $
and $\tilde{g}_{\ogamma}(y_1, \ldots, y_{r_1+r_2} ) = (x_1, \ldots x_{n_{\K}}) $.
We now observe that
 ${\tilde{g}_{\ogamma}}^{-1}$ is continuous on $\R^{r+1} \times [0, 2\pi)^{r_2}$. 
In fact, $ [\frac{1}{2},1] \times \prod_{i=1}^{r} [0,\frac{1}{m_i}] \times (0,2\pi)^{r_2}$
is homeomorphic to its image in $\R_{\overline{\gamma}}^{r_1}  \times  \C^{r_2}$
under the map ${\tilde{g}_{\ogamma}}^{-1}$.
Therefore 
\begin{equation*}
\d{\tilde{g}_{\ogamma}}^{-1} \left( [\frac{1}{2},1] 
\times \prod_{i=1}^{r} [0,\frac{1}{m_i}] \times (0,2\pi)^{r_2} \right)
 =
{\tilde{g}_{\ogamma}}^{-1} \d \left( [\frac{1}{2},1] 
\times \prod_{i=1}^{r} [0,\frac{1}{m_i}] \times (0,2\pi)^{r_2} \right).
\end{equation*}
It is easy to see that
\begin{equation*}
\d{\tilde{g}_{\ogamma}}^{-1} \left( [\frac{1}{2},1] \times \prod_{i=1}^{r} [0,\frac{1}{m_i}] \times [0,2\pi)^{r_2} \right)
 ~\subseteq~~
{\tilde{g}_{\ogamma}}^{-1} \d \left( [\frac{1}{2},1] \times \prod_{i=1}^{r} [0,\frac{1}{m_i}] \times [0,2\pi)^{r_2} \right).
\end{equation*}
We now define $2(r+1)$ sets  as follows; $I_{1,1} = \{\frac{1}{2} \}~\times~\prod_{i=1}^{r}~[0, \frac{1}{m_i}]$,
$I_{1,2} =  \{1\}~\times~\prod_{i=1}^{r}~[0,\frac{1}{m_i}]$ and for $2 \le j \le r+1$
\begin{equation*}
I_{j, l} 
= 
\begin{cases} 
[\frac{1}{2},1]~\times~\prod_{k=1}^{j-2}~[0,\frac{1}{m_i}] 
\times \{0\} \times ~\prod_{k=j}^{r}~[0,\frac{1}{m_i}] 
& \text{ if    } l=1 \\ 
[\frac{1}{2},1]~\times~\prod_{k=1}^{j-2}~[0,\frac{1}{m_i}] \times \{\frac{1}{m_{j-1}}\} 
\times ~\prod_{k=j}^{r}~[0,\frac{1}{m_i}]
& \text{ if    } l=2.
\end{cases}
\end{equation*}
For $1\le j \le r+1$, define $\psi_{j,1},  \psi_{j,2} :  [0,1]^{n_{\K}-1} \to  I_{j, 1} \times [0,2\pi]^{r_2}$
as follows; 
\begin{eqnarray*}
\psi_{j,1} (t_1, \cdots, t_{n_{\K}-1})
&=&
\begin{cases} 
(\frac{1}{2},  \frac{t_1}{m_1},  \cdots  \frac{t_{r}}{m_r}, 2\pi t_{r+1} 
\cdots, 2\pi t_{n_{\K}-1}) & \text{ if   } j=1, \\
(\frac{1 + t_1}{2},  \frac{t_2}{m_1}, \cdots \frac{t_{j -1}}{m_{j - 2}}, 
 0, \frac{t_{j}}{m_j} \cdots  \frac{t_{r}}{m_r}, 2\pi t_{r+1} \cdots, 2\pi t_{n_{\K}-1})
&\text{  if    } j > 1,
\end{cases} \\
\text{and } 
\psi_{j,2} (t_1, \cdots, t_{n_{\K}-1})
&=&
\begin{cases} 
(1,  \frac{t_1}{m_1},  \cdots  \frac{t_{r}}{m_r}, 2\pi t_{r+1} \cdots, 2\pi t_{n_{\K}-1}) 
& \text{ if  } j = 1, \\
(\frac{1 + t_1}{2},  \frac{t_2}{m_1}, \cdots \frac{t_{j-1}}{m_{j-2}}, 
\frac{1}{m_{j-1}}, \frac{t_{j}}{m_j} \cdots  \frac{t_{r}}{m_r}, 2\pi t_{r+1} \cdots, 2\pi t_{n_{\K}-1}) 
& \text{ if  } j > 1.
\end{cases}
\end{eqnarray*}
Further, for $\ogamma \in \{ \pm 1 \}^{r_1}$, we define 
$h_{\ogamma} :  I_{j, 1} \times [0,2\pi]^{r_2} \to \R_{\ogamma}^{r_1} \times \C^{r_2}$
and $h'_{\ogamma} :  \R_{\ogamma}^{r_1} \times \C^{r_2} \to \R^{n_{\K}}$ 
as follows; $h_{\ogamma}(x_1, \cdots, x_{n_{\K}}) = (y_1, \cdots, y_{r_1+ r_2})$,
where $y_i$'s for $1 \le i  \le r_1$ are defined by \eqref{inv1} and $y_i$'s for 
$r_1 + 1 \le i  \le r_1+ r_2$  are defined by \eqref{inv2}
and $h'_{\ogamma}(y_1, \cdots,  y_{r_1+r_2}) = (z_1, \cdots, z_{n_{\K}} )$, where
for $1 \le i \le r_1$, $z_i = y_i$ and for $r_1+1 \le i \le r_1 + r_2$, 
$$
z_i = |y_i|\cos (\arg(y_i))
\phantom{m}\text{and}\phantom{m}
z_{r_2+i} =  |y_i| \sin (\arg(y_i)),
$$
where argument for $y_i$'s are inside $[0, 2\pi)$.  We now define $2(r+1)$ maps
$\phi_{j, l}  : [0,1]^{n_{\K} -1} \to \R^{n_{\K}}$ for $1\le j \le r+1$ and
$l= 1, 2$ as follows; 
$\phi_{j, l} = h'_{\ogamma} \circ {h}_{\overline{\gamma}} \circ \psi_{j, l}$.
These $\phi_{j, l}$ cover the boundary of 
$(h'_{\ogamma} \circ \tilde{g}_{\overline{\gamma}}^{-1} )\left( [\frac{1}{2},1] \times 
\prod_{i=1}^{r} [0,\frac{1}{m_i}] \times [0,2\pi)^{r_2} \right)$ as $h'_{\ogamma}$
is a homeomorphism from $\R_{\ogamma}^{r_1} \times \C^{r_2}$ onto its image
inside $\R^{n_{\K}}$ and hence
$$
\d \left(h'_{\ogamma} \left( \tilde{g}_{\overline{\gamma}}^{-1}\left( [\frac{1}{2},1] \times 
\prod_{i=1}^{r} [0,\frac{1}{m_i}] \times [0,2\pi)^{r_2} \right) \right)\right)
=
h'_{\ogamma} \left(\d \tilde{g}_{\overline{\gamma}}^{-1} \left( [\frac{1}{2},1] \times 
\prod_{i=1}^{r} [0,\frac{1}{m_i}] \times [0,2\pi)^{r_2} \right) \right).
$$
We now compute the constants in \lemref{lipschitz} for each of the maps $\phi_{j, l}$. 
Since we know that 
$m_{k-1} = \max_i \ceil{\log |\sigma_i(\eta_{k-1})|}$ and $0 \le t_{k-1} \le 1$, we have
$$
 |\sigma_i(\eta_{k-1})|^{\frac{t_{k-1}}{m_{k-1}}}  ~\le~   e 
 $$ 
for $1 \le i \le r_1$ and $2 \le k \le r+1$.
Further for any $t, t' \in [0, 1]$, using mean value theorem, we get
$$
\left| |\sigma_i(\eta_{k-1})|^{\frac{t}{m_{k-1}}}
 - |\sigma_i(\eta_{k-1})|^{\frac{t'}{m_{k-1}}} \right| 
 ~\le~ 
 e ~ |t - t'|.
$$
From now on, we denote the $i$-th projection map of $\phi_{j, l}$ 
by $\phi^{i}_{j, l}$. Applying \lemref{lipschitz} and the above
observations, for any $\overline{t},  \overline{t}' \in [0,1]^{n_{\K} -1}$
and $1 \le i \le r_1$, we have
$$
\left| \phi^{i}_{1, 1}( \overline{t})
 - \phi^{i}_{1, 1} ( \overline{t}' ) \right| 
 ~\le~ 
 \left(\frac{1}{2} \right)^{\frac{1}{n_{\K}}} r e^r ~ || \overline{t} - \overline{t}'  ||.
$$
Since for any $t, t' \in [0,1]$
\begin{equation}\label{lip-bd}
|\cos(2\pi t) - \cos(2\pi t')| ~\le~  2\pi ~ |t - t'|
\phantom{m}\text{and}\phantom{m}
|\sin(2\pi t) - \sin(2\pi t')| ~\le~  2\pi ~ |t - t'|,
\end{equation}
therefore for any $r_1 +1 \le i \le n_{\K}$ and for any 
$\overline{t},  \overline{t}' \in [0,1]^{n_{\K} -1}$,  as before using \lemref{lipschitz},
we have
$$
\left| \phi^{i}_{1, 1}( \overline{t})
 - \phi^{i}_{1, 1} ( \overline{t}' ) \right| 
 ~\le~ 
 \left(\frac{1}{2} \right)^{\frac{1}{n_{\K} }} (2\pi + r)~ e^r 
 ~ || \overline{t} - \overline{t}'  ||.
$$
Now applying second part of \lemref{lipschitz}, for any 
$\overline{t},  \overline{t}' \in [0,1]^{n_{\K} -1}$, we have
\begin{equation}\label{lip-1.1}
\left| \phi_{1, 1}( \overline{t})
 - \phi _{1, 1} ( \overline{t}' ) \right| 
 ~\le~ 
 \left(\frac{1}{2} \right)^{\frac{1}{n_{\K}}} \sqrt{n_{\K}}~(2\pi + r)~ e^r 
 ~ || \overline{t} - \overline{t}'  ||.
\end{equation}
Proceeding similarly, for any $\overline{t},  \overline{t}' \in [0,1]^{n_{\K} -1}$, we get
\begin{equation}\label{Lip-1.2}
\left| \phi_{1, 2}( \overline{t})
 - \phi _{1, 2} ( \overline{t}' ) \right| 
 \le 
\sqrt{n_{\K}}~ (2\pi + r)~ e^r 
~  || \overline{t} - \overline{t}'  ||.
\end{equation}
Since for any $t,  t' \in [0,1]$
\begin{equation}\label{lip-bd2}
\left|\left(\frac{1 + t}{2}\right)^{\frac{1}{ n_{\K}  }}
-
 \left(\frac{1 + t'}{2}\right)^{\frac{1}{ n_{\K}  }} \right| 
~\le~
\frac{2^{\frac{-1}{n_{\K} }}}{ n_{\K} }~ |t - t'|,
\end{equation}
we have for any $2 \le j \le 2(r+1)$, $1\le l\le 2$, $\overline{t},  \overline{t}' \in [0,1]^{n_{\K} -1}$
and $1 \le i \le r_1$
$$
\left| \phi^{i}_{j, l}( \overline{t})
 - \phi^{i}_{j, l} ( \overline{t}' ) \right| 
 \le
 \left( \frac{2^{\frac{-1}{ n_{\K} }}}{n_{\K}} + r - 1 \right) e^{r}
~ || \overline{t} - \overline{t}'  ||.
$$
Using \eqref{lip-bd} and \eqref{lip-bd2} and proceeding as before,
we get for any $2 \le j \le 2(r+1)$, $1\le l\le 2$, 
$\overline{t},  \overline{t}' \in [0,1]^{n_{\K} -1}$ 
and $r_1 + 1 \le i \le n_{\K}$
$$
\left| \phi^{i}_{j, l}( \overline{t})
 - \phi^{i}_{j, l} ( \overline{t}' ) \right| 
 ~\le~ 
\left( 2\pi + \frac{2^{\frac{-1}{n_{\K}}}}{n_{\K}} + r - 1 \right)
e^{r}
|| \overline{t} - \overline{t}'  ||
~\le~
(2\pi +r) e^r  ~|| \overline{t} - \overline{t}'  ||.
$$
Combining, for any $2 \le j \le 2(r+1)$, $1\le l\le 2$ and 
$\overline{t},  \overline{t}' \in [0,1]^{n_{\K} -1}$, we have 
\begin{equation}\label{Lip-j.l}
\left| \phi_{j, l}( \overline{t})
 - \phi _{j, l} ( \overline{t}' ) \right| 
 \le 
\sqrt{n_{\K}}~ (2\pi + r) e^r 
  ~|| \overline{t} - \overline{t}'  ||.
\end{equation}
 If the interval $[\frac{1}{2},1)$ is replaced by $[\frac{t^{n_{\K}}}{2}, t^{n_{\K}})$ for some 
positive $t \in \R$, we deduce in a similar way that the bound $L$ in the definition \ref{def-2}
of Lipschitz class is less than $t\sqrt{n_{\K}} (2\pi+r) e^r$.
Hence $\d\F_{1/2, \overline{\gamma}}(0,\frac{1}{m_1}, \cdots, 0 \frac{1}{m_r}, t^{n_{\K}})$
is of Lipschitz class $\L(n_{\K}, 2r+2, Lt)$, where $L = \sqrt{n_{\K}} (2\pi + r) e^r$.
\end{proof}

\subsection{Preliminary lemmas}
Using \thmref{KZ} and \thmref{Dob}, we derive the following bounds
on the product $\prod_{j=1}^r m_j$, where $m_{j} = \max_i \ceil{\log |\sigma_i(\eta_{j})|}$.
\begin{lem}\label{Regbound}
 Let $\q_1$ be a modulus of $\K$. There exist $r$ units $\eta_1, \cdots \eta_r$ 
modulo roots of unity which generate $U_{\q_1}$ such that
 $$
\frac{R_{\K, \q_1}}{2^r (r+1)^{\frac{r-1}{2}}} 
~\le~ 
\prod_{j=1}^{r} m_j  
~\le~ 
 7^r (r+1)^{r+ 1/2} n_{\K}^{2r} ~R_{\K,\q_1},
 $$
 where as before $r= r_1 + r_2 -1$ and $R_{\K,\q_1}$ is the $\q_1$ regulator of the field $\K$.
\end{lem}
 
\begin{proof}
If there are no fundamental units, the regulator $R_{\K,\q_1}$ is defined
to be $1$ and hence the above inequalities are trivially satisfied.
Hence from now on, we assume that $r \ge 1$.
Let $\{ \eta'_1, \cdots \eta'_r\}$ be any set of generators of $U_{\q_1}$
modulo roots of unity.  Let $l$ denote the map
\begin{eqnarray*}
l : \R^{r_1} \times \C^{r_2} & \to & \R^{r+1} \\
(x_i)_i & \to & (e_i \log |x_i|)_{i=1}^{r + 1}
\end{eqnarray*}
Here $e_i$ is $2$ if $r_1+ 1 \le i \le r+1$
and $1$ otherwise.  Set $\Lambda_r =  l (\phi(U_{\q_1}))$, where 
$\phi$ is as in subsection 5.1. Note that $\Lambda_r$ is a 
lattice of rank $r$. The vectors
$$
\vec{v_1} 
=
 \left(\frac{1}{\sqrt{r+1}}, \cdots ,\frac{1}{\sqrt{r+1}}\right), 
 \phantom{m}
\vec{v_j} 
=  
\left(e_i\log |\sigma_i(\eta'_{j-1})|\right)_{i=1}^{r+1}
$$
for $2 \le j \le r +1$ form a basis for $\R^{r+1}$.
Since $||\vec{v_1} || =1$ and $v_1$ is orthogonal to the
$\R$ vector space generated by $\{ v_2, \cdots v_{r+1} \}$, the volume of the
lattice generated by the vectors $\{ v_1, \cdots v_{r+1} \}$ in $\R^{r+1}$
is the same as the volume of $\Lambda_r$ in $\R^r$.
Hence the volume of $\Lambda_r$ is $\sqrt{r+1} ~R_{\K,{\q_1}}$.
Note that the volume of $\Lambda_r$ is independent 
of the choice of the generators $\{ \eta'_1, \cdots \eta'_r\}$
of $U_{\q_1}$.

Let $\{\vec{w_1} \ldots \vec{w_{r}}\}$ be a reduced Korkin Zolotarev
basis for $\Lambda_r$. Choose a set of generators $\eta_1, \cdots \eta_r$ of $U_{\q_1}$ 
modulo roots of unity such that $f( \phi(\eta_i) ) = \vec{w_i}$. Then for $1 \le j \le r$
$$
m_j 
~\ge~
 \frac{1}{2} \text{max}_i( e_i\log(|\sigma_i(\eta_j)|) ) 
~\ge~ 
\frac{1}{2\sqrt{r+1}}||\vec{w_j}||.
$$
Since $\prod_{j=1}^r ||\vec{w_j}|| $ is greater than or equal to the volume, we have
$$
\prod_{j=1}^{r} m_j 
~\ge~   
\left(\frac{1}{2\sqrt{r+1}} \right)^r \text{ Vol }(\Lambda_r) 
 ~\ge~
 \frac{R_{\K, \q_1}}{2^r (r+1)^{\frac{r-1}{2}}} .
 $$
We now compute the upper bound. Since there is at least one
fundamental unit, we known that $n_{\K} > 1$. If $n_{\K}=2$, 
then $R_{\K,\q_1} = |\log |\eta|~|$, where
$\eta$ generates $U_{\q_1}$ modulo roots of unity.
Since $R_{\K,\q_1} \ge R_{\K} \ge 1/5$ (see \cite{Friedman*89}),
we get the required upper bound in this case.
If $n_{\K} \ge 3$, then
$$
1 + \frac{\log n_{\K}}{6n_{\K}^2} 
~\ge~ 
1 + \frac{1}{7n_{\K}^2-1} 
~\ge~ 
e^{\frac{1}{7n_{\K}^2}}.
$$
Now by Dobrowolski's theorem, we know that $\max_{i}( e_i\log |\sigma_i(\eta_j)| )
\ge \frac{1}{7n_{\K}^2}$. This implies that $m_j ~\le~  7n_{\K}^2 \max_i(e_i \log|\sigma_i(\eta_j)|) 
 ~\le~ 7n_{\K}^2 ~||\vec{w_j}||$. Using \thmref{KZ}, we now get
$$
\prod_{j=1}^r m_j 
~\le~ (7n_{\K}^2)^r \prod_{j=1}^r ||\vec{w_j}|| 
~\le~
(7n_{\K}^2)^r\left(\frac{3+r}{4}\right)^{r/2} r^{r/2}~\sqrt{r+1}~ R_{\K,\q_1}
~\le~  7^r (r+1)^{r+ 1/2} n_{\K}^{2r} ~R_{\K,\q_1}.
$$
\end{proof}

\begin{lem}\label{IN}
Let $\q_1$ be a modulus of $\K$ with $n_{\K} > 1$, $\{\sigma_1, \cdots, \sigma_{r} \}$ be a set of embeddings
of $\K$ into $\C$ such that no two embeddings are conjugate to each other. 
For $\eta_j$ and $m_j$ as defined earlier and $\alpha \in \K$, the integral
\begin{equation}\label{int-1}
\int_{\R^r} \frac{dx_1\cdots dx_r}{\text{max}_{1 \le i \le r+1}(\prod_{j=1}^r |\sigma_i(\eta_j)|^{x_j/m_j}|\sigma_i(\alpha)|)^{n_{\K}-1}} 
= 
 \frac{m_1\cdots m_r}{R_{\K,\q_1}|\N(\alpha)|^{\frac{n_{\K}-1}{n_{\K}}}} 
 \left( \frac{n_{\K}}{n_{\K}-1} \right)^r.
\end{equation}
\end{lem}

\begin{proof}
Using the transformation $ x_j \to m_j x_j$, we see that the left hand side of \eqref{int-1} is
equal to
\begin{equation}\label{int-2}
m_1 \cdots m_r ~\int_{\R^r} \frac{dx_1\cdots dx_r}{\text{max}_{1 \le i \le r+1}(\prod_{j=1}^r 
|\sigma_i(\eta_j)|^{x_j}|\sigma_i(\alpha)|)^{n_{\K} - 1}}.
\end{equation}
We now make the substitution $y_i = \sum_{j=1}^r x_j \log|\sigma_i(\eta_j)|$,
for $1 \le i \le r$. Hence \eqref{int-2} is equal to
$$
m_1 \cdots m_r ~\frac{\prod_{i=1}^r e_i}{R_{\K,\q_1}}
\int_{\R^r} \frac{dy_1\cdots dy_r}{\text{max}( e^{y_1}|\sigma_1(\alpha)|, \cdots e^{y_r}
|\sigma_r(\alpha)|, e^{-Y} |\sigma_{r+1}(\alpha)| )^{n_{\K}-1}}
$$
where $e_i$ for $1 \le i \le r$ is $1$ or $2$ depending on $\sigma_i$ is real or 
complex embedding and
\begin{align*}
Y 
& =  - \sum_{j=1}^r x_j \log |\sigma_{r+1}(\eta_j)| 
 = \frac{1}{e_{r+1}} \sum_{j=1}^r  x_j  (-e_{r+1} \log |\sigma_{r+1}(\eta_j)| )\\
& =  \frac{1}{e_{r+1}} \sum_{j=1}^r x_j \sum_{k=1}^r e_k \log |\sigma_k(\eta_j)|
 =   \frac{1}{e_{r+1}} \sum_{k=1}^r e_k y_k.
\end{align*}
The integral is now identical to the one computed in the proof
of Lemma 10 of \cite{Debaene*19}.
\end{proof}

\subsection{Counting points in the fundamental domain}

Applying \thmref{Widmer}, we will now derive the following
counting theorem.
\begin{thm}\label{counting}
Let $\a, \q$ be co-prime ideals of $\O_{\K}$, $\ogamma \in \{\pm 1 \}^{r_1}$ and
$\mathfrak{C}$ be the ideal class of $\a\q$ in the class group of $\O_{\K}$.
For $\vec{k} \in \prod_{j=1}^{r} ([0,m_j) \cap \Z)$, let 
$\beta_{\vec{k}} = \left( \prod_{j=1}^r |\sigma_i(\eta_j)|^{-k_j/m_j}\right)_{i=1}^{r+1}$
and $\Lambda_{n_{\K}}(\vec{k})$ be the lattice $h' \left(\phi(\a\q) \cdot \beta_{\vec{k}}\right)$
in $\R^{n_{\K}}$, where $\phi, h'$ are as in subsection 5.3.  Also let
$$
S\left(\a, \q,  \F_{\frac{1}{2}, \ogamma}(t^{n_{\K}}) \right)  
= \{\alpha \in \a~|~\phi(\alpha) \in \F_{\frac{1}{2}, \ogamma}(t^{n_{\K}}), ~\alpha \equiv b \bmod \q\},
$$
where $b \in \O_{\K}$.
Then for any real number $t \ge 1$, we have
\begin{equation*}
\left| S\left(\a, \q,  \F_{\frac{1}{2}, \ogamma}(t^{n_{\K}}) \right)  \right| 
 = 
\frac{(2 \pi)^{r_2} R_{\K,\q_1} t^{n_{\K}} }{\sqrt{4|d_{\K}|}~ \N(\a\q)} ~ + ~
\rO^*\left( \frac{ e^{n_{\K}^2 + 8 n_{\K}} n_{\K}^{ \frac{3}{2}n_{\K}^2  + \frac{11}{2}n_{\K} -\frac{1}{2}}t^{n_{\K}-1} }
{\mathfrak{N}(\mathfrak{C}^{-1})^{-1} |\N(\a\q)|^{\frac{n_{\K}-1}{n_{\K}}}} + m_1\cdots m_r\right).
\end{equation*}
where $\q_1 = \q \q_{\infty}$ with $\q_{\infty}$ containing all the infinite places of $\K$ and  
$$
\N(\mathfrak{C}^{-1}) = 
\text{max}_{\b_i \in \mathfrak{C}^{-1}} 
\sum_{i =1}^{m_1\cdots m_r} \frac{1}{|\mathfrak{N}(\b_i)|^{\frac{n_{\K}-1}{n_{\K}}}}.
$$
The term $m_1\cdots m_r$ may be omitted when
$
t ~\ge~ \frac{ \max_{\vec{k}}~\delta_1\left( \Lambda_{n_{\K}} (\vec{k}) \right) }{\sqrt{n_{\K}}(2\pi +r)e^r}$
or $\q = \O_{\K}$.
\end{thm}

\begin{proof}
Since we want to count $\alpha \in \a$ such that $\alpha \equiv b \bmod \q$, where $\a, \q$
are co-prime, we need to count $\alpha$ in only one residue class, say $a$ modulo $\a\q$.
Recall that $\K \neq \Q$.  We need to count 
$$
\left|(\phi(\a\q) + \phi(a)) \bigcap \F_{\frac{1}{2},\overline{\gamma}}(t^{n_{\K}}) \right|.
$$
By \lemref{SPT}, we have
\begin{equation}\label{partition}
| ( \phi(\a\q) + \phi(a) ) \bigcap \F_{{\frac{1}{2},\ogamma}}(t^{n_{\K}})| 
= 
\sum_{\vec{k}} | \left((\phi(\a\q) + \phi(a)) \cdot \beta_{\vec{k}} \right) \bigcap \F_{ \frac{1}{2}, \ogamma }
(0, \frac{1}{m_1},\cdots, 0, \frac{1}{m_r}, t^{n_{\K}} ) |,
\end{equation}
where ${\vec{k}}  \in \prod_{j=1}^{r} ([0,m_j) \cap \Z)$. Since 
$$
\text{Vol}(h' \left( (\phi(\a\q)+ \phi(a)) \cdot \beta_{\vec{k}} \right) ) 
~=~
\text{Vol}(h'(\phi(\a\q)\cdot \beta_{\vec{k}})),
$$
applying Widmer's estimate (\thmref{Widmer}), the main term for each $\vec{k}$ is 
$$
\frac{\text{Vol} \left( h'\left( \F_{\frac{1}{2},\ogamma}(0,\frac{1}{m_1},
\cdots, 0, \frac{1}{m_r}, t^{n_{\K}} )\right) \right)}
{\text{Vol}\left( h' \left( \phi(\a\q)\cdot \beta_{\vec{k}} \right) \right)}.
$$
The volume of $h'(\phi(\a\q))$ is the determinant of the matrix whose column
vectors form a basis for this lattice $h'(\phi(\a\q))$. If $\beta_{\vec{k}} = (\beta_i)_{i=1}^{r+1}$,
multiplying the $i$-the entry of each basis vector with $\beta_i$ for $1 \le i \le r_1+r_2$
and for $r_1 + r_2 \le i \le r_1 + 2r_2$ with $\beta_{ i - r_2}$, we get that
$$
\text{Vol}(h'(\phi(\a\q)\cdot \beta_{\vec{k}})) = \text{Vol}(h'(\phi(\a\q))).
$$
This implies that the main term for each $\vec{k}$ is independent of $\vec{k}$ and
equal to
\begin{equation}\label{vol-1}
\frac{\text{Vol}(h'(\F_{\frac{1}{2}, \ogamma }(0,\frac{1}{m_1},\cdots, 0, \frac{1}{m_r}, t^{n_{\K}})))}
{\text{Vol}(h'(\phi(\a\q)))}.
\end{equation}
We know that $\text{Vol}(h'(\phi(\a\q))) = 2^{-r_2} \sqrt{|d_{\K}|} \N(\a\q)$ (see page 31 of \cite{Neukirch*99}).
To compute the volume
$\text{Vol}(h'(\F_{\frac{1}{2}}(0,\frac{1}{m_1},\cdots, 0, \frac{1}{m_r}, t^{n_{\K}})))$,
we note that 
\begin{equation*}
\text{Vol} \left(h' \left(\F_{\frac{1}{2}} 
(0,\frac{1}{m_1},\cdots, 0, \frac{1}{m_r}, ~t^{n_{\K}} )\right)\right) 
 = 
\int_{h'(\F_{\frac{1}{2}} (0,\frac{1}{m_1}, \cdots, 0, \frac{1}{m_r}, ~t^{n_{\K}}))} 
dx_1 \cdots dx_{r_1} dx_{r_1+1} \cdots dx_{n_{\K}},
\end{equation*}
where the variables $x_{r_1+1}$ to $x_{n_{\K}}$ are complex. By definition
$$
\F_{\frac{1}{2}} \left(0,\frac{1}{m_1},\cdots, 0, \frac{1}{m_r}, t^{n_{\K}} \right) 
=  
g^{-1}\left((\frac{t^{n_{\K}}}{2}, t^{n_{\K}}] \times \prod_{j=1}^r [0, \frac{1}{m_j})\right), 
$$
the argument for each complex coordinate in the
pre-image covers the entire interval $[0,2\pi)$.
Replacing the complex variables
with polar co-ordinates and integrating over the
arguments,
\begin{eqnarray*} 
&& 
\text{Vol} \left(h' \left(\F_{\frac{1}{2}}(0,\frac{1}{m_1},\cdots, \frac{1}{m_r}, t^{n_{\K}} )\right) \right) \\
& = &
2^{r_1} (2\pi)^{r_2} \int_{|x_1|} d|x_1| \cdots  \int_{|x_{r_1}|} 
~d|x_{r_1}|  \int_{|x_{r_1 +1}|} |x_{r_1+1}| ~~d|x_{r_1 +1}| \cdots \int_{|x_{r+1}|}  |x_{r + 1}| ~~d|x_{r+1}|.
\end{eqnarray*}
The ranges in the above integral are clear from the formulae \eqref{id} with
$\alpha(x) \in (\frac{t^{n_{\K}}}{2}, t^{n_{\K}}]$ and $\alpha_j(x) \in [0, \frac{1}{m_j})$. 
We now make a change of variable $x= (x_1,  x_2, \cdots x_{r+1})$ by 
$(\alpha(x), \alpha_1(x), \cdots \alpha_{r}(x))$.
To compute the Jacobian, we note that for $1 \le i \le r+1$ and $1 \le j \le r$, we have
\begin{equation*}
\frac{\partial |x_i|}{\partial \alpha(x)} = \frac{|x_i|}{n_{\K}\alpha(x)}
\phantom{m}\text{ and }\phantom{m}
\frac{\partial |x_i|}{\partial \alpha_j(x)} = |x_i| \log |\sigma_i(\eta_j)|. 
\end{equation*}
We also note that $\alpha(x) = \prod_{i=1}^{r+1} |x_i|^{e_i}$
where $e_i=2$ for complex coordinate and $e_i =1$ otherwise.
Now rewriting the integral, we have
\begin{align*} 
\text{Vol} \left(\F_{\frac{1}{2}} (0,\frac{1}{m_1},\cdots, 0, \frac{1}{m_r}, t^{n_{\K}} )\right) 
& = 
2^{r_1} (2\pi)^{r_2} 2^{-r_2} R_{\K,\q_1} \int_{0}^{1/m_1} d\alpha_1(x) \cdots 
\int_{0}^{1/m_r} d\alpha_r(x) 
\int_{ \frac{t^{n_{\K}} }{2}  }^{ t^{n_{\K} } } d\alpha(x)  \\
& = 
 \frac{2^{r_1-1} \pi^{r_2} R_{\K,\q_1} t^{n_{\K}} }{m_1 \cdots m_r}.
\end{align*}
Therefore for each $\vec{k}$, the main term \eqref{vol-1} is equal to
$$
\frac{1}{2^{r_1}} \cdot \frac{2^{r_1-1} \pi^{r_2} R_{\K,\q_1} t^{n_{\K}} }{m_1 \cdots m_r} 
\cdot \frac{1}{2^{-r_2} \sqrt{|d_{\K}|} \N(\a\q)} 
= 
\frac{(2\pi)^{r_2} R_{\K,\q_1} t^{n_{\K}}}{ (m_1 \cdots m_r) \sqrt{4|d_{\K}|} \N(\a\q)}.
$$
Hence the main term after summing over $\vec{k} \in \prod_{j=1}^r [0,m_j) \cap \Z$ is
equal to
\begin{equation}\label{vol-main}
\frac{(2 \pi)^{r_2} R_{\K,\q_1}t^{n_{\K}} }{\sqrt{ 4 |d_{\K}|} \N(\a\q)} .
\end{equation}
Applying \thmref{Widmer} and \lemref{Lip-const}, the error term for each $\vec{k}$
is bounded by
\begin{equation}\label{vol-er}
(2r+ 2) n_{\K}^{3n_{\K}^2/2} \text{max}_{0 \le i < n_{\K}} 
\frac{(\sqrt{n_{\K}}(2\pi + r)e^r t)^i}{\prod_{j=1}^i  \delta_j( \Lambda_{n_{\K}} (\vec{k}) )}
\le 
(2r+ 2) n_{\K}^{3n_{\K}^2/2} \text{max}_{0 \le i < n_{\K}} 
\left(\frac{\sqrt{n_{\K}}(2\pi + r)e^r t}{ \delta_1( \Lambda_{n_{\K}} (\vec{k}) )  }\right)^i.
\end{equation}
If we have  $\sqrt{n_{\K}}(2\pi + r)e^r t \ge \delta_1( \Lambda_{n_{\K}} (\vec{k}) )$, then 
to get an upper bound of \eqref{vol-er}, we can replace $i$ by $n_{\K}-1$.

To deduce the asymptotic for all $t \ge 1$, we first consider the case where
$\delta_1( \Lambda_{n_{\K}} (\vec{k}))~>~\sqrt{n_{\K}}(2\pi + r) e^rt$
for  some vector $\vec{k}$. In this case, we claim that 
$$
 \left| \left( (\phi(\a\q) + \phi(a)) \cdot \beta_{\vec{k}} \right) \bigcap \F_{ \frac{1}{2}, \ogamma }
\left(0, \frac{1}{m_1},\cdots, 0, \frac{1}{m_r}, t^{n_{\K}} \right) \right| \le 1.
$$
Suppose not and let $x,y \in \left((\phi(\a\q) + \phi(a)) \cdot \beta_{\vec{k}} \right) 
\bigcap \F_{ \frac{1}{2}, \ogamma } \left(0, \frac{1}{m_1},\cdots, 0, \frac{1}{m_r}, t^{n_{\K}} \right)$ 
be distinct. Then $x-y \in \phi(\a\q) \cdot \beta_{\vec{k}}$. By definition, this implies that 
$\delta_1( \Lambda_{n_{\K}} (\vec{k}) ) \le ||h'(x-y)||$. Applying \lemref{length}, 
we get
$$
\sqrt{n_{\K}}(2\pi + r)~e^r t
~<~
 \delta_1( \Lambda_{n_{\K}} (\vec{k}) ) 
 ~\le~
|| h' (x-y) || 
~\le~
2 (\sqrt{r+1}~) ~e^r t,
$$
a contradiction. If $\q = \O_{\K}$, the same argument applies to $x$
in place of $x-y$ and this implies that there are no exceptional points.
The total number of $\vec{k}$ for which $ \delta_1( \Lambda_{n_{\K}} (\vec{k}) )
~>~\sqrt{n_{\K}}(2\pi + r)e^r t$ is at most $m_1\cdots m_r$.
Hence
$$
\sum_{\vec{k} \atop{\sqrt{n_{\K}}(2\pi + r)e^r t ~<~\delta_1( \Lambda_{n_{\K}} (\vec{k}) )}} 
 \left| \left((\phi(\a\q) + \phi(a)) \cdot \beta_{\vec{k}} \right) \bigcap \F_{ \frac{1}{2}, \ogamma }
\left(0, \frac{1}{m_1},\cdots, 0, \frac{1}{m_r}, t^{n_{\K}} \right) \right| 
~\le~ 
m_1\cdots m_r.
$$
For each $\vec{k}$ for which $\sqrt{n_{\K}}(2\pi + r)e^r t  < \delta_1( \Lambda_{n_{\K}} (\vec{k}) )$, 
we claim that 
$$
\frac{(2\pi)^{r_2} R_{\K,\q_1} t^{n_{\K}}}{ (m_1 \cdots m_r) \sqrt{4|d_{\K}|} \N(\a\q)}
~\le~ 
(2r+ 2) n_{\K}^{3n_{\K}^2/2} 
\left(\frac{\sqrt{n_{\K}}(2\pi + r)e^r t}{ \delta_1( \Lambda_{n_{\K}} (\vec{k}) )  }\right)^{n_{\K} -1}.
$$
Note that $2^{r-1} \pi^{r_2} (r + 1)^{\frac{r-1}{2}} n_{\K}^{n_{\K}}
~\le~
(2r + 2) n_{\K}^{\frac{n_{\K}}{2} (3n_{\K} + 1) } (2\pi + r)^{n_{\K}}$
when $n_{\K} \ge 2$ and hence
$$
\frac{(2\pi)^{r_2} 2^{r} (r+1)^{\frac{r-1}{2}} t^{n_{\K}}}{2\sqrt{|d_{\K}|} \N(\a\q)}
 ~\le~ 
 \frac{(2r+2) n_{\K}^{\frac{3n_{\K}^2}{2}} }{   \text{Vol}(h'(\phi(\a\q)\cdot \beta_{\vec{k}})) } 
 \left(  \frac{ \sqrt{n_{\K}} (2\pi + r) t }{ n_{\K} } \right)^{n_{\K}}.
$$
By \thmref{KZ}, we have $\delta_1(\Lambda_{n_{\K}}(\vec{k}))^{n_{\K}} ~\le~ 
n_{\K}^{n_{\K}} \text{Vol}(h'(\phi(\a\q)\cdot \beta_{\vec{k}}))$
and therefore
\begin{equation*}
\frac{(2\pi)^{r_2} 2^r (r+1)^{\frac{r-1}{2}} t^{n_{\K}}}{2 \sqrt{|d_{\K}|} \N(\a\q)}
 \le (2r+2)n_{\K}^{\frac{3n_{\K}^2}{2} } 
 \left( \frac{\sqrt{n_{\K}}(2\pi + r) e^r t}{\delta_{1}(\Lambda_{n_{\K}}(\vec{k}))}\right)^{n_{\K}}.
\end{equation*}
Finally applying \lemref{Regbound} and $\sqrt{n_{\K}}(2\pi + r)e^r t 
< \delta_1( \Lambda_{n_{\K}} (\vec{k}) )$, we get
$$
\frac{(2\pi)^{r_2} R_{\K,\q_1} t^{n_{\K}}}{m_1\cdots m_r\sqrt{4|d_{\K}|} \N(\a\q)}
 \le (2r+2)n_{\K}^{\frac{3n_{\K}^2}{2} } 
 \left( \frac{\sqrt{n_{\K}}(2\pi + r) e^r t}{\delta_{1}(\Lambda_{n_{\K}}(\vec{k}))}\right)^{n_{\K}-1},
$$
as claimed. Therefore for $t \ge 1$ the error is bounded by
\begin{equation}\label{Error1}
(2r+ 2) n_{\K}^{\frac{3n_{\K}^2}{2}} (\sqrt{n_{\K}}(2\pi + r)e^r t)^{n_{\K}-1}
\sum_{\vec{k}} \frac{1}{ \delta_1( \Lambda_{n_{\K}} (\vec{k}) )^{n_{\K}-1}  }
~+~~  m_1\cdots m_r.
\end{equation}
If there are no $\vec{k}$ such that $\sqrt{n_{\K}}(2\pi +r)e^rt ~<~ 
\text{ max}_{\vec{k}}~\delta_1\left( \Lambda_{n_{\K}} (\vec{k}) \right)$, 
then the error term is bounded by
\begin{equation}\label{Error}
(2r+ 2) n_{\K}^{3n_{\K}^2/2} (\sqrt{n_{\K}}(2\pi + r)e^r t)^{n_{\K}-1}
\sum_{\vec{k}} \frac{1}{ \delta_1( \Lambda_{n_{\K}} (\vec{k}) )^{n_{\K}-1}  }.
\end{equation}
Let 
$$
\mu(\a\q, \vec{k}) = \text{min}_{\alpha \in \a\q} \text{max}_{1\le i \le r+1} 
\left( |\sigma_i(\alpha)| \prod_{j=1}^r |\sigma_i(\eta_j)|^{-k_j/m_j} \right).
$$
From the definition of successive minima $\delta_1( \Lambda_{n_{\K}} (\vec{k}) ) \ge \mu(\a\q, \vec{k})$.
For any $\alpha \in \a\q$, let ${\rm K}_{\alpha}$ be the set of all $\vec{k}$ 
for which the minimum $\mu(\a\q, \vec{k})$ is attained at $\alpha$, i.e.
$$
{\rm K}_{\alpha}= \left\{ \vec{k} ~~\Big| ~  \mu(\a\q, \vec{k}) 
= \text{max}_{1 \le i \le r+1} (|\sigma_i(\alpha)| 
\prod_{j=1}^r |\sigma_i(\eta_j)|^{-k_j/m_j}) \right\}.
$$
We also set $Y_1 = \{ \alpha \in \a\q ~|~ {\rm K}_{\alpha} \ne \emptyset \}$.
Substituting in \eqref{Error} (analogously in \eqref{Error1}), we get
\begin{align*}
\sum_{\vec{k}} \frac{1}{\delta_1( \Lambda_{n_{\K}} (\vec{k}) )^{n_{\K}-1}} 
&\le  
\sum_{\vec{k}} \frac{1}{\mu(\a\q, \vec{k})^{n_{\K} -1}}
=
\sum_{\alpha \in Y_1} \sum_{\vec{k} \in {\rm K}_{\alpha}}  \frac{1}{\mu(\a\q, \vec{k})^{n_{\K} -1}}\\
& \le  
\sum_{\alpha \in Y_1}' \sum_{u \in U_{\q_1}} \sum_{\vec{k} \in {\rm K}_{u\alpha}} 
\frac{1}{\text{max}_{1 \le i \le r+1}\left(|\sigma_i(u\alpha)| \prod_{j=1}^r 
|\sigma_i(\eta_j)|^{-k_j/m_j}\right)^{n_{\K}-1}} \\
& \le  
\sum_{\alpha \in Y_1}' \sum_{\vec{k} \in \Z^r} \frac{1}{\text{max}_{1 \le i\le r+1} 
(|\sigma_i(\alpha)| \prod_{j=1}^r |\sigma_i(\eta_j)|^{-k_j/m_j})^{n_{\K}-1}},
\end{align*}
where $'$ on the inner sum indicates that the sum is over non-associate elements
of $\a\q$ with respect to the unit group $U_{\q_1}$. Note that there are at most $m_1\cdots m_r$
elements in the outermost sum.  Applying \lemref{IN}, we now bound the inner sum to get
\begin{equation*}
\sum_{\vec{k} \in \Z^r} \frac{1}{\text{max}_{1\le i \le r+1} (|\sigma_i(\alpha)| 
\prod_{j=1}^r |\sigma_i(\eta_j)|^{-k_j/m_j})^{n_{\K}-1}}
~\le~
\frac{2^{r} m_1\cdots m_r}{ R_{\K,\q_1} |\N(\alpha)|^{\frac{n_{\K} -1}{n_{\K} }}} .
\end{equation*}
Now the term
\begin{equation*}
\sum_{\alpha \in Y_1 }'  \frac{1}{|\mathfrak{N}(\alpha)|^{\frac{n_{\K}-1}{n_{\K} }}} 
\le 
\text{max}_{\alpha _i \in \a\q} \sum_{i =1}^{m_1\cdots m_r} 
 \frac{1}{|\mathfrak{N}(\alpha_i)|^{\frac{n_{\K}-1}{n_{\K}}}} 
\le
\left(\frac{1}{\N(\a\q)}\right)^{\frac{n_{\K} -1}{n_{\K}}} 
\text{max}_{\b_i \in \mathfrak{C}^{-1}} 
\sum_{i =1}^{m_1\cdots m_r} \frac{1}{|\mathfrak{N}(\b_i)|^{\frac{n_{\K}-1}{n_{\K}}}} .
\end{equation*}
The inner sum is a constant that depends only on the
inverse of the class of $\a\q$ in the class group. We shall denote this 
by $\mathfrak{N}(\mathfrak{C}^{-1})$. Combining all these estimates and
applying \lemref{Regbound}, we get that the error term \eqref{Error} is bounded by
\begin{align*}
& 
2^{r+1} (r+ 1)^{2r +1}  (10n_{\K})^{2r}
n_{\K}^{\frac{3n_{\K}^2}{2}} (\sqrt{n_{\K}}(2\pi + r)e^r)^{n_{\K}-1}    
\frac{\mathfrak{N}(\mathfrak{C}^{-1}) t^{n_{\K}-1} }{|\N(\a\q)|^{\frac{n_{\K}-1}{n_{\K}}}} \\
& \le~
 e^{n_{\K}^2 + 8 n_{\K}} n_{\K}^{ \frac{3}{2}n_{\K}^2  + \frac{11}{2}n_{\K} -\frac{1}{2}}
 \frac{\mathfrak{N}(\mathfrak{C}^{-1}) t^{n_{\K}-1} }{\N(\a\q)^{\frac{n_{\K}-1}{n_{\K}}}}.
\end{align*}
\end{proof}
\smallskip
\begin{rmk}
One can replace the condition $t ~\ge~ \frac{ \max_{\vec{k}}~\delta_1\left( \Lambda_{n_{\K}}
(\vec{k}) \right) }{\sqrt{n_{\K}}(2\pi +r)e^r}$ in Theorem~\ref{counting}  by
\begin{equation*}
  t\ge \frac{\left(\left(\frac{2}{\pi}\right)^{r_2} \sqrt{|d_{\K}|} ~\N(\a\q)\right)^{\frac{1}{n_{\K}}} e^{\sum_{j=1}^r (m_j - 1)}}{ \sqrt{n_{\K}(r+1)}}.
\end{equation*}


Indeed, by Minkowski's lattice point theorem, there exists a point 
$a \in \a\q$ such that
$$
\phi(a)\cdot \beta_{\vec{k}} \in \phi(\a\q) \cdot \beta_{\vec{k}} 
\phantom{m}\text{and}\phantom{m}
\N(a) <  \left(\frac{2}{\pi} \right)^{r_2} \sqrt{|d_{\K}|} ~\N(\a\q).
$$
By \eqref{id}, we have
$$
|\sigma_i (a)| = \N(a)^{\frac{1}{n_{\K}}} \prod_{j=2}^{r+1} |\sigma_i(\eta_{j-1})|^{\alpha_j(a)}.
$$
Let $b_j = \floor{\alpha_j(a) - \frac{k_{j-1}}{m_{j-1}}}$ and $\tilde{b_j} = \alpha_j(a) - \frac{k_{j-1}}{m_{j-1}} - b_j$.
We observe that $ 0 \le \tilde{b_j} < 1$.
Consider the unit $u_{a} = \prod_{j=2}^{r+1} \eta_{j-1}^{b_{j}}$ in $\O_{\K}^*$.
We now have
\begin{eqnarray*}
|| \phi(a \cdot u_{a}^{-1}) \cdot \beta_{\vec{k}} ||
& = &
\sqrt{ \sum_{i=1}^{r+1} |\sigma_i( a \cdot u_{a}^{-1}) |^2 \prod_{j=2}^{r+1}
 |\sigma_i(\eta_{j-1})|^{-\frac{2k_{j-1}}{m_{j-1}}}} \\
& = &
\N(a)^{\frac{1}{n_{\K}}} 
\sqrt{ \sum_{i=1}^{r+1} \prod_{j=2}^{r+1} |\sigma_i(\eta_{j-1})|^{2 \alpha_j(a) - 2b_j -\frac{2k_{j-1}}{m_{j-1}}}} \\
& = &
\N(a)^{\frac{1}{n_{\K}}} 
\sqrt{ \sum_{i=1}^{r+1} \prod_{j=2}^{r+1} |\sigma_i(\eta_{j-1})|^{2\tilde{b_j}}} \\
& \le &
\N(a)^{\frac{1}{n_{\K}}}  \sqrt{r+1}~ e^{\sum_{j=1}^r m_j}. 
\end{eqnarray*}
\end{rmk}

\subsection{Counting ideals in ray classes}
In this subsection, we complete the proof of  \thmref{asymfinal}. 
We use notations from the previous sections. We start by
simplifying the main term in Theorem~\ref{counting}.

\begin{lem}{\rm (Debaene \cite{Debaene*19},  Lemma 12)}\label{bd-De}
Let $\b_1, \b_2, \cdots$ be integral ideals in $\O_{\K}$, ordered such that  $\N(\b_1) \le \N(\b_2)  \cdots$.
Then for any real number $y \ge 2$
$$
\sum_{i=1}^y \N(\b_i)^{ \frac{1}{n_{\K}}  -1} 
~\le~ 
6n_{\K} y^{\frac{1}{n_{\K}}} (\log y)^{\frac{(n_{\K}-1)^2}{n_{\K}}}.
$$
\end{lem}
\begin{lem}[Lang \cite{Lang*94}, page 127]\label{classnumber}
Let $\q_1 = \q\q_{\infty}$ be a modulus of $\K$, $r_1, U_{\q_1}, h_{\K}$ be as defined earlier 
and $h_{\K, \q}$ be the cardinality of the narrow ray class group of $\K$ modulo $\q$.  
Then
$$
h_{\K, \q} = \frac{ 2^{r_1}  \varphi(\q) h_{\K} }{ [\O_{\K}^*  :  U_{\q_1}] } ~.
$$
\end{lem}
\begin{lem}
  \label{mainterm}
 Let $r_1, r_2, h_{\K},  d_{\K}, R_\K,  \q_1, U_{\q_1}, h_{\K, \q}, R_{\K, \q_1}$ be as
 in the previous subsections. Also let $\alpha_{\K}$ be the residue of the Dedekind zeta function
 at $s=1$, $\mu_{\K}$ be the group of
 roots of unity in $\K$ and $\mu_{\q_1}$ be the cardinality of $U_{\q_1} \cap \mu_{\K}$.
 We have
  \begin{equation*}
    \frac{(2 \pi)^{r_2} R_{\K, \q_1}}{\mu_{\q_1}\sqrt{|d_{\K}|}}
    =
    \frac{\alpha_{\K} \varphi(\q)}{h_{\K,\q}},
    \quad\quad
    \frac{R_{\K,\q_1}}{R_\K}=\frac{\mu_{\q_1}}{|\mu_\K|}
    \frac{2^{r_1}\varphi(\q)h_\K}{h_{\K,\q}}.
  \end{equation*}
\end{lem}
\begin{proof}
Applying the analytic class number formula \eqref{acf}, we see that  
\begin{equation*}
\frac{(2 \pi)^{r_2} R_{\K,\q_1}}{\mu_{\q_1}\sqrt{|d_{\K}|} \N(\q)} 
= 
\frac{2^{r_1}(2 \pi)^{r_2} h_{\K} R_{\K}}{\mu_{\q_1}\sqrt{|d_{\K}|} \N(\q)} 
 \times \frac{R_{\K,\q_1}}{ 2^{r_1} h_{\K} R_{\K}}
=
\frac{\alpha_{\K}}{\N(\q)} \times \frac{|\mu_{\K}| R_{\K,\q_1}}{ 2^{r_1}  \mu_{\q_1}  h_{\K}  R_{\K} }.
\end{equation*}
As shown in \lemref{Regbound}, consider the lattice $ l( \phi(U_{\q_1} ))$ and 
note that $\sqrt{r+1} ~R_{\K,\q_1}$ is the volume of its fundamental domain.
Also note that $\sqrt{r+1} ~R_{\K}$ is the volume of the fundamental domain
of $l(\phi(\O_{\K}^*))$.  We know that $l(\phi(U_{\q_1}))$ and $l(\phi(\O_{\K}^*))$
are finitely generated free modules over $\Z$ of the same rank. Therefore by the
structure theorem for finitely generated free modules over principal ideal domains, 
we have
$$
 [l(\phi(\O_{\K}^*)) : l(\phi(U_{\q_1}))]
= \frac{\text{Vol}( l(\phi(U_{\q_1}))}{\text{Vol}(l(\phi(\O_{\K}^*)))} 
= \frac{R_{\K,\q_1}}{R_{\K}}.
$$
Along with the torsion part coming from the 
roots of unity, we have
$$
[\O_{\K}^* : U_{\q_1}]  = [l(\phi(\O_{\K}^*)) : l(\phi(U_{\q_1}))] \frac{|\mu_{\K}|}{\mu_{\q_1}}.
$$
Applying the above identities, we have
$$
\frac{\alpha_{\K}}{\N(\q)} \times \frac{|\mu_{\K}| R_{\K,\q_1}}{ 2^{r_1}  \mu_{\q_1}  h_{\K}  R_{\K} }
=
\frac{\alpha_{\K}}{\N(\q)} \times\frac{ [\O_{\K}^*  :  U_{\q_1}] }{ 2^{r_1} h_{\K}  }   ~.
$$
Now using \lemref{classnumber}, we get the first formula.
The second one follows along similar lines;
\begin{equation*}
  \frac{R_{\K,\q_1}}{R_\K}
  =
  \frac{\mu_{\q_1}}{|\mu_\K|}[\O_{\K}^*  :  U_{\q_1}]
  =
  \frac{\mu_{\q_1}}{|\mu_\K|}
  \frac{2^{r_1}\varphi(\q)h_\K}{h_{\K,\q}}.
\end{equation*}
This completes the proof of this lemma.
\end{proof}

We can now proceed to the proof of Theorem~\ref{asymfinal}.
\begin{proof}[Proof of Theorem~\ref{asymfinal}]
  Let us fix an ideal $\c \in [\b]^{-1}$. Since $\c\b = (\alpha)$ for
  some $\alpha \equiv 1 \bmod^* \q$, in order to count the number of
 integral ideals in $[\b]$ with norm at most $x$, it is sufficient to count
  $(\alpha), \alpha \in \c$ such that $\alpha \equiv 1 \bmod^* \q$ of
  norm at most $x\N\c$. This implies that
\begin{equation*}
  \sum_{\substack{ \a\subset\O_\K,\\ [\a]=[\b] \\ \N\a \le x} }1
  =
  \frac{1}{\mu_{\q_1}} \left| 
    \left\{ \alpha \in \c  ~|~ \phi(\alpha) \in
      \F_{\overline{\eta}}(0,1, \cdots, 0,1,x\mathfrak{N}\c), 
      ~\alpha \equiv 1 \bmod \q \right\} 
  \right|,
\end{equation*}
where $\overline\eta = (1, \cdots, 1)$. Let $\Lambda_{n_{\K}}(\vec{k})$ be the lattice
$h' \left(\phi(\c\q) \cdot \beta_{\vec{k}}\right)$ in $\R^{n_{\K}}$.
Using \thmref{counting}, we know that if
$x \ge 1$, we get
\begin{equation}\label{eqasym}
  \sum_{\substack{\a\subset\O_\K,\\ [\a]=[\b] \\ \N\a \le x}}1
  = \frac{(2 \pi)^{r_2} R_{\K, \q_1}}{\mu_{\q_1}\sqrt{|d_{\K}|}~ \N\q} x
  +
  \rO^* \left( \frac{e^{n_{\K}^2 + 8 n_{\K}} }{\mu_{\q_1}}
    n_{\K}^{ \frac{3}{2}n_{\K}^2  + \frac{11}{2}n_{\K} -\frac{1}{2} } 
    \N(\mathfrak{C}^{-1}) 
    \left(\frac{x}{\N\q}\right)^{1-\frac{1}{n_\K}}
    + \frac{m_1\cdots m_r}{\mu_{\q_1}} \right).
\end{equation}
We rewrite the main term by appealing to the first part of
Lemma~\ref{mainterm}.  We then bound above $m_1\cdots m_r/\mu_{\q_1}$
first by the quantity $ 7^r  n_{\K}^{4n_{\K}} R_{\K,\q_1}/\mu_{\q_1}$ thanks to
Lemma~\ref{Regbound}. We bound it further by
$ 2^{3n_{\K}} n_{\K}^{4n_{\K}}2^{r_1}\varphi(\q)\frac{R_\K}{|\mu_\K|}\frac{h_\K}{h_{\K,\q}}$ by
invoking the second part of Lemma~\ref{mainterm}. The paper
\cite{Friedman*89} by E.~Friedman ensures us that
$R_{\K}/|\mu_\K| \ge 0.2$, and so, this upper bound is at least equal
to~2. Lemma~\ref{bd-De} may thus be applied to majorize
$\N(\mathfrak{C}^{-1})$.  
On shortening $F(\q)$ by $F$, then $n_\K$ by $n$ and $R_\K/|\mu_\K|$
by $\tilde{R}$, this leads to the error term
\begin{equation*}
   e^{n^2 + 8 n} 
   n^{ \frac{3}{2}n^2  + \frac{11}{2}n -\frac{1}{2} }
   6n \bigl((2n)^{4n}F\tilde{R} \bigr)^{1/n}
   \log\bigl((2n)^{4n}F\tilde{R} \bigr)^{\frac{(n-1)^2}{n}}
   \left(\frac{x}{\N\q}\right)^{1-\frac{1}{n}}
   ~+~ (2n)^{4n}F\tilde{R} .
\end{equation*}
We separate the contribution of $\q$ and of $\K$ by using
$$
\log\bigl((2n)^{4n}F\tilde{R} \bigr)\le
2\log\bigl((2n)^{4n}\tilde{R} \bigr)\log(3F)
$$
which is a consequence of the inequality $a+b\le 2ab$ valid when
$a,b\ge1$.
But we first notice that $\log((2n)^{4n}F\tilde{R})\ge 1$ so that we
may simply replace the exponent $(n-1)^2/n$ by $n$.
We further check that
\begin{equation*}
  e^{n^2 + 8 n} 
  n^{ \frac{3}{2}n^2  + \frac{11}{2}n -\frac{1}{2} }
  6n(2n)^4
  \le 500\,n^{12n^2}.
\end{equation*}
Hence the theorem follows.
\end{proof}

\bigskip
\noindent
{\bf Acknowledgements.} 
Research of this article was partially supported by Indo-French 
Program in Mathematics (IFPM).  All authors would like to thank IFPM 
for financial support. The first author would also like to 
acknowledge MTR/2018/000201, SPARC project~445 and DAE number 
theory plan project for partial financial support.

We also thank Ethan Lee from Canberra for pointing out
to us the work \cite{Sunley*73} of J.~Sunley.

\bigskip
\bibliographystyle{plain}
\bibliography{\jobname.bib}

\begin{thebibliography}{10}

\bibitem{Debaene*19}
K.~Debaene.
\newblock Explicit counting of ideals and a {B}run-{T}itchmarsh inequality for
  the {C}hebotarev density theorem.
\newblock {\em Int. J. Number Theory}, 15(5):883--905, 2019.

\bibitem{Dobrowolski*79}
E.~Dobrowolski.
\newblock On a question of {L}ehmer and the number of irreducible factors of a
  polynomial.
\newblock {\em Acta Arith.}, 34(4):391--401, 1979.

\bibitem{Friedman*89}
E.~Friedman.
\newblock Analytic formulas for the regulator of a number field.
\newblock {\em Invent. Math.}, 98(3):599--622, 1989.

\bibitem{Lagarias-Lenstra-Schnorr*90}
J.~C. Lagarias, H.~W. Lenstra~Jr., and C.-P. Schnorr.
\newblock Korkin-{Z}olotarev bases and successive minima of a lattice and its
  reciprocal lattice.
\newblock {\em Combinatorica}, 10(4):333--348, 1990.

\bibitem{Lang*94}
S.~Lang.
\newblock {\em Algebraic number theory}, volume 110 of {\em Graduate Texts in
  Mathematics}.
\newblock Springer-Verlag, New York, second edition, 1994.

\bibitem{LEE22}
E.~S. Lee.
\newblock On the number of integral ideals in a number field.
\newblock {\em Journal of Mathematical Analysis and Applications},
  517(1):126585, 2023.

\bibitem{Neukirch*99}
J.~Neukirch.
\newblock {\em Algebraic number theory}, volume 322 of {\em Grundlehren der
  mathematischen Wissenschaften [Fundamental Principles of Mathematical
  Sciences]}.
\newblock Springer-Verlag, Berlin, 1999.
\newblock Translated from the 1992 German original and with a note by Norbert
  Schappacher, With a foreword by G. Harder.

\bibitem{Sunley*71}
J.~E.~S. Sunley.
\newblock {\em On the class numbers of totally imaginary quadratic extensions
  of totally real fields}.
\newblock ProQuest LLC, Ann Arbor, MI, 1971.
\newblock Thesis (Ph.D.)--University of Maryland, College Park.

\bibitem{Sunley*73}
J.~E.~S. Sunley.
\newblock Class numbers of totally imaginary quadratic extensions of totally
  real fields.
\newblock {\em Trans. Amer. Math. Soc.}, 175:209--232, 1973.

\bibitem{Sutherland*15}
A.~Sutherland.
\newblock Class field theory, ray class groups and ray class fields.
\newblock MIT Mathematics, 18.785, Number Theory I, Lecture \# 20, 2015.

\bibitem{Tatuzawa*73}
T.~Tatuzawa.
\newblock On the number of integral ideals in algebraic number fields, whose
  norms not exceeding {$x$}.
\newblock {\em Sci. Papers College Gen. Ed. Univ. Tokyo}, 23:73--86, 1973.

\bibitem{Widmer*10a}
M.~Widmer.
\newblock Counting primitive points of bounded height.
\newblock {\em Trans. Amer. Math. Soc.}, 362(9):4793--4829, 2010.

\end{thebibliography}

\end{document}